\documentclass{amsart}
\usepackage{graphicx,subfigure}
\usepackage[english]{babel}
\usepackage[cp1251]{inputenc}
\usepackage{amssymb}
\usepackage{lscape}
\usepackage{amsthm}
\usepackage{inputenc}

\newtheorem{defn}{Definition}[section]
\newtheorem{thm}{Theorem}[section]
\newtheorem{pr}{Proposition}[section]
\newtheorem{exam}{Example}[section]
\newtheorem{cor}{Corollary}[section]

\newenvironment{dem}{\rm \trivlist \item[\hskip \labelsep{\it
      Proof}:]}{\par\nopagebreak \hfill $\Box$ \endtrivlist}

\def\ll{\mathcal{L}}
\def\CC{\mathbb{C}}

\begin{document}
\date{}
\author{L.M. Camacho, E.M. Ca\~{n}ete, J.R. G\'{o}mez, B.A. Omirov}

\title[3-filiform Leibniz algebras of maximum length]{\bf 3-filiform Leibniz algebras of maximum length.}
\thanks{L.M. Camacho and E.M. Ca\~nete would like to thank to
the Institute of Ma\-the\-ma\-tics and Information Technologies of Uzbekistan for their hospitality. The second author was partially supported by a Grant of the Junta de Andaluc\'{i}a (Spain).}

\maketitle
\begin{abstract}
This work completes the study of the solvable Leibniz algebras, more precisely, it completes the classification of the $3$-filiform Leibniz algebras of maximum length \cite{3-filiform}. Moreover, due to the good structure of the algebras of maximum length, we also tackle some of their cohomological properties.
Our main tools are the previous result of Cabezas and Pastor \cite{Pastor}, the construction of appropriate homogeneous basis in the considered connected gradation and the computational support provided by the two programs implemented in the software \textit{Mathematica}.
\end{abstract}



\section{Introduction}
Leibniz algebras appear from the cohomology study done by Loday in 1993 \cite{loday1} and they are further investigated by several authors as
Ayupov, Casas and others (\cite{Ayupov}, \cite{Casas}). In the cohomology study there is an important family of Leibniz algebras: those whose length of the gradation is maximum. The remarkable fact that an algebra can be decomposed into direct sum of subspaces of dimension 1 makes easier the calculations of the derivations since they induce the corresponding gradation of the group of cohomologies.

The main goal of this paper is to continue the study of the $p$-filiform Leibniz algebras of maximum length. These algebras play a main role in mathematics over the last years, either in the classification theory or in geometrical, analytical and physical applications.

  In early works we have already closed the classification of the $p$-filiform Leibniz algebras of maximum length for $0 \leq p \leq 2$ (see \cite{Ayupov}, \cite{J.Lie.Theory2}). Here we study the 3-filiform Leibniz algebras of maximum length, their spaces of derivations and their first cohomology group.

Moreover, we will use three programs very helpful to obtain the classification of maximum length algebras and their space of derivations.

 Recall \cite{loday1} that an algebra $\ll$ over a field $F$ is called a Leibniz algebra if it satisfies the following Leibniz identity:
$$[x,[y,z]]=[[x,y],z]-[[x,z],y], \forall x,y,z \in \ll $$
where $[.,.]$ denotes the multiplication in $\ll$.

Consider an arbitrary algebra $\ll$ in the set of n-dimensional Leibniz algebras over a field $F$.
Let $B=\{ e_1,$ $ e_2,$ $\cdots$ $e_n\}$ be a basis of $\ll$. Then $\ll$ is determined, un to isomorphisms, by the multiplication
rule for the basis elements; namely,
$$[e_i,e_j]=\displaystyle\sum_{k=1}^n \gamma_{ij}^k e_k$$
where $\gamma_{ij}^k$ are the structure constants. Therefore, fixing a basis, we can regard each algebra of dimension
n over a field F as a point in the $n^3$-dimensional space of structure constants endowed with the Zariski
topology.

 From now on the Leibniz algebras will be considered over the field of complex numbers $\mathbb{C}$, and with finite dimension. Let $\ll$ be a Leibniz algebra, then $\ll$ is naturally filtered by the descending central sequence $\ll^1=\ll$, $\ll^{k+1}=[\ll^k,\ll]$ with $k\geq 1.$ Thus, a nilpotent algebra $\ll$ has nilindex equal to $s$ if $s$ is the minimum integer such that $\ll^s \neq \{0\}$ and $\ll^{s+1}=\{0\}.$

 A Leibniz algebra $\ll$ is $\mathbb{Z}$-graded if $\ll=\oplus_{i \in \mathbb{Z}}V_i$,
   where $[V_i,V_j]\subseteq V_{i+j}$ for any $i,j \in \mathbb{Z}$ with a finite number of non null spaces $V_i$.

  We will say that a $\mathbb{Z}$-graded Leibniz algebra $\ll$ admits a \emph{connected gradation} if $\ll=V_{k_1}\oplus V_{k_1+1} \oplus \dots \oplus V_{k_1+t}$ and $V_{k_1+i}\neq <0>$ for any $i$ $(0 \leq i \leq t)$.

 Let us define the naturally graded algebras as follows:
   \begin{defn}
   Let us take $\ll_i=\ll^i/\ll^{i+1}$, $1\leq i \leq k$ and $gr \ll=\ll_1 \oplus \ll_2 \oplus \dots \oplus \ll_k$. Then $[\ll_i,\ll_j]\subseteq \ll_{i+j}$
    and we obtain the graded algebra $gr \ll$. If $gr \ll$ and $\ll$ are isomorphic, in notation $gr\ll \cong \ll$, we say that $\ll$
    is a naturally graded algebra.
   \end{defn}

    The above constructed gradation is called \emph{natural gradation}.

    \begin{defn}\label{def:length}
  The number $l( \oplus \ll)=l(V_{k_1}\oplus V_{k_1+1} \oplus \dots \oplus V_{k_1+t})=t+1$ is called the
length of the gradation, where $ \oplus \ll$ is a connected gradation. The gradation $\oplus \ll$ has maximum length if $l(\oplus \ll)=dim( \ll)$.
  \end{defn}

  We define the length of an algebra $\ll$ by:
  \begin{center}
 $l(\ll)=\max \{l(\oplus \ll)  \hbox{ such that } \oplus \ll = V_{k_1}\oplus \dots \oplus V_{k_t}\hbox{ is a connected gradation}\}.$
 \end{center}

An algebra $\ll$ is called of maximum length if $l(\ll)=dim (\ll)$.


 The set $R(\ll)=\{x \in \ll: [y,x]=0, \ \forall y \in \ll\}$ is called \emph{the right annihilator of $\ll$}. $R_x$ denotes the operator $R_x: \ll \rightarrow \ll$ such that $R_x(y)=[y,x], \ \forall y \in \ll$ and it is called the right operator. The set $Cent(\ll)=\{z \in \ll: [x,z]=[z,x]=0, \ \forall x \in \ll\}$ is called \emph{the center of $\ll$}.

  Let $x$ be a nilpotent element of the set $\ll \setminus \ll^2$. For the nilpotent operator $R_x$ we define a descending sequence $C(x)=(n_1,n_2, \dots, n_k)$, which consists of the dimensions of the Jordan blocks of the operator $R_x$. In the set of such sequences we consider the lexicographic order, that is,
  $C(x)=(n_1,n_2, \dots, n_k)< C(y)=(m_1, m_2, \dots, m_s)$ if and only if there exists $i \in \mathbb{N}$ such that $n_j=m_j$ for any $j<i$ and $n_i<m_i$.

  \begin{defn}\label{def:char.seq}
  The sequence $C(\ll)=\max C(x)_{x \in \ll \setminus \ll^2}$ is called the characteristic sequence of the algebra $\ll$.
  \end{defn}

  Let $\ll$ be an $n$-dimensional nilpotent Leibniz algebra and $p$ a non negative integer ($p<n$).
  \begin{defn}\label{def:p-fili}
  The Leibniz algebra $\ll$ is called $p$-filiform if $C(\ll)=(n-p,\underbrace{1,\dots,1}_{p})$. If $p=0$, $\ll$ is called null-filiform
  and if $p=1$ it is called filiform.
  \end{defn}
  Therefore, an algebra with the characteristic sequence $(n-2,1,1)$ is called 2-filiform, whereas a nilpotent algebra with nilindex $n-2$ is called quasi-filiform. Note that in the Lie algebras case both definitions coincide.

   \begin{defn}\label{def:derivation}
A linear transformation $d$ of a Leibniz algebra $\ll$ is called a derivation of $\ll$ if
$$d([x,y])=[d(x),y]+[x,d(y)] \hbox{ for any } x,y \in \ll.$$
Denote by $Der(\ll)$ the set of all derivations.
  \end{defn}

  It is clear that the right operator $R_x $ is a derivation for any $x \in \ll$. Derivations of this type are called inner derivations. Similar to the Lie algebras case the set of the inner derivations forms an ideal of the algebra $Der(\ll).$

  Since our algebra is $\mathbb{Z}-$graded, i.e $\ll=\oplus_{i \in \mathbb{Z}} V_i$, this gradation induces a gradation of the algebra $Der(\ll)=\oplus_{i \in \mathbb{Z}} W_i$ in the following way:
$$W_i=\{d_i \in Der(\ll): d_i(x) \in V_{i+j} \hbox{ for any } x \in V_j \}.$$

\

For an $n$-dimensional algebra of maximum length it is easy to see that $Der(\ll)=W_{-n} \oplus \dots \oplus W_{n}$ (see \cite{Omirov2}).
For more details see the definition of the cohomology groups for Leibniz algebras introduced in \cite{loday-piras}.

\section{$3$-filiform non-Lie Leibniz algebras of maximum length}

 In this section we are going to continue the classification of the $p$-filiform Leibniz algebras of ma\-xi\-mum length. The study of the filiform and 2-filiform cases has been already done in \cite{J.Lie.Theory2}, so we are going to continue with the 3-filiform Leibniz algebras case.

The used technique in this section is as follows: we will extend the naturally graded 3-filiform Leibniz algebras by using the natural gradations. In this way, we can distinguish two cases: the natural graded Lie algebras and the natural graded non-Lie algebras. The study of the first case was closed in \cite{3-filiform}, so we explain the results obtained in the second family. After that we will work with a homogeneous basis and we will assume that the associated gradation has maximum length. Finally, we will use some programs implemented in the software \textit{Mathematica} (which will be explained below) as well as properties of the gradation and of the nilpotence to arrive at a contradiction or at the classification.

I would like to stress in the fact that using computer programs is very helpful to achieve the presented classification. Two programs will be used in this section: the program of the Leibniz identity and the program of isomorphisms. The first program computes the Leibniz identity of a Leibniz algebra and was presented in \cite{JSC}. The second one establishes when two algebras are isomorphic, moreover we have added some subroutines to know if two algebras are isomorphic or not, when one of them is an uniparameter family. It returns the value of the parameter for which would be isomorphics. The algorithmic method can be found on \cite{mitesis}.

 The implementation of these programs are presented in low and fixed dimension. Then we will formulate the generalizations, proving by induction the results for arbitrary fixed dimension. Finally, point out that the algorithmic method of these programs are presented with a step-by-step explanation in the following Web site: http://personal.us.es/jrgomez.

\

\subsection{Non split case}

\

In this section we will restrict our study to classify the 3-filiform Leibniz algebras of maximum length, which are the extension of the non split and naturally graded non-Lie Leibniz algebras. The Lie case has been closed in \cite{3-filiform}, where there is not any non split 3-filiform Leibniz algebra of maximum length.

First of all, let us see the classifications of the naturally graded 3-filiform Leibniz algebras (\cite{NGp-F}).

\begin{thm}
Let $\ll$ be a complex $n$-dimensional non-split naturally graded $3$-filiform non-Lie Leibniz algebra and $n \geq 7$. Then there exists a basis  $\{e_1,e_2,\dots, e_{n-3},f_1,f_2,f_3\}$ of the algebra, such that $\ll$ is isomorphic to
$$L^{1}:\begin{cases}
                [e_{i},e_{1}]=e_{i+1},  &1\leq i \leq n-4, \\
                [e_{1},f_{1}]=f_{3}, \\
                [e_{i},f_{2}]=e_{i+1},  &1\leq i \leq n-4.
          \end{cases}$$
\end{thm}

\begin{thm}
Let $\ll$ be a complex $n$-dimensional non split $3$-filiform non-Lie Leibniz algebra and $n \geq 7$. Then $l(\ll)<n.$
\end{thm}

\begin{dem}

The natural gradation of $L^1$ is: $\ll_1\oplus \dots\oplus \ll_{n-3}$ where $\ll_1=<e_1,f_1,f_2>$, $\ll_2=<e_2,f_3>$ and $\ll_i=<e_i>$ for $3\leq i \leq n-3$. We are going to study the length of its extension, which is denoted by $\widetilde{L^1}.$ Point out that we call \textit{the extension of the algebra} as the natural generalization of the structural constants of the algebra, using the information of its associated natural gradation.

 Note that $\{e_2,e_3, \dots,e_{n-3},f_3\}$ belong to the ideal $R(L^1)$ and $e_{n-3} \in Cent(L^1)$. Moreover $[e_1,f_1]+[f_1,e_1] \in R(L^1)$, then we conclude $f_3 \in R(L^1)$. Finally, by taking the change of basis $e_1^{'}=e_1$,  $e_{i+1}^{'}=[e_i^{'},e_1^{'}]$ for $1 \leq i \leq n-4$, $f_1^{'}=f_1$ and $f_2^{'}=f_2$, we can write the law of $\widetilde{L^1}$ as:
$$\begin{cases}
[e_i,e_1]=e_{i+1}, &1\leq i \leq n-4,\\
[e_1,f_1]=f_3+(*)e_3+\dots+(*)e_{n-3},\\
[e_i,f_2]=e_{i+1}+(*)e_{i+2}\dots+(*)e_{n-3}, &1\leq i \leq n-4,\\
[f_i,e_1]=(*)e_3+\dots+(*)e_{n-3}, &1\leq i\leq 2,\\
[f_3,e_1]=(*)e_4+\dots+(*)e_{n-3},\\
[f_i,f_j]=(*)e_3+\dots +(*)e_{n-3}, & 1\leq i,j \leq 2,\\
[f_3,f_i]=(*)e_4+\dots+(*)e_{n-3}, & 1 \leq i \leq 2\\
[e_i,f_1]=(*)e_{i+2}+ \dots+ (*)e_{n-3}, & 2 \leq i \leq n-5,
\end{cases}$$
where the asterisks $(*)$ denote the corresponding coefficients in the products.
A crucial tool in the proof is the construction of a homogeneous basis, which generators are:
\begin{align*}
\widetilde{x_s}&=e_1+ \sum_{i=2}^{n-3}a_ie_i +\sum_{j=1}^{3}a_{n-3+j}f_j,\\
\widetilde{x_t}&=f_1+ \sum_{i=1}^{n-3}b_ie_i +\sum_{j=2}^{3}b_{n-3+j}f_j,\\
\widetilde{x_u}&=f_2+ \sum_{i=1}^{n-3}c_ie_i +\sum_{j=1, j\neq 2}^{3}c_{n-3+j}f_j.
\end{align*}

Therefore the products of the generators of $\widetilde{L^1}$ can be defined in the new basis as follows:

\begin{align*}
[\widetilde{x_s},\widetilde{x_s}]&=(1+a_{n-1})e_2+(*)e_3+ \dots + (*)e_{n-3}+a_{n-2}f_3,\\
[\widetilde{x_t},\widetilde{x_t}]&=b_1(b_1+b_{n-1})e_2+(*)e_3 + \dots+ (*)e_{n-3}+b_1f_3,\\
[\widetilde{x_u},\widetilde{x_u}]&=c_1(1+c_1)e_2+ (*)e_3+ \dots+ (*)e_{n-3}+c_1 c_{n-2} f_3,\\
[\widetilde{x_s},\widetilde{x_t}]&=(b_1+b_{n-1})e_2+(*)e_3 + \dots+ (*)e_{n-3}+f_3,\\
[\widetilde{x_t},\widetilde{x_s}]&=b_1(1+a_{n-1})e_2+(*)e_3 + \dots+ (*)e_{n-3}+b_1a_{n-2}f_3, \\
[\widetilde{x_s},\widetilde{x_u}]&=(1+c_1)e_2+(*)e_3 + \dots+ (*)e_{n-3}+c_{n-2}f_3, \\
[\widetilde{x_u},\widetilde{x_s}]&= c_1(1+a_{n-1})e_2+(*)e_3 + \dots+ (*)e_{n-3}+c_1a_{n-2}f_3,\\
[\widetilde{x_t},\widetilde{x_u}]&=b_1(1+c_1)e_2+(*)e_3 + \dots+ (*)e_{n-3}+b_1c_{n-2}f_3, \\
[\widetilde{x_u},\widetilde{x_t}]&=c_1(b_1+b_{n-1})e_2+(*)e_3 + \dots+ (*)e_{n-3}+c_1f_3.\\
\end{align*}

Since $\{\widetilde{x_s},\widetilde{x_t},\widetilde{x_u}\}$ are linearly independent, then $det\left(\begin{array}{ccc}
1 & a_{n-2} & a_{n-1}\\
b_1 & 1 & b_{n-1}\\
c_1 & c_{n-2} & 1
\end{array}\right)\neq 0.$

\


$\bullet$ Case 1: If $1+a_{n-1} \neq 0,$ we have the following subcases:

\

\fbox{Case 1.1:} If $[\widetilde{x_s},\widetilde{x_s}]$ and $[\widetilde{x_s},\widetilde{x_t}]$ are linearly independent, we take the homogeneous basis $y_1=\widetilde{x_s},$ $y_i=[y_{i-1},y_1]$ for $2 \leq i \leq n-3,$ $z_1=\widetilde{x_t},$ $z_2=\widetilde{x_u}$ and $z_3=[\widetilde{x_s},\widetilde{x_t}]=[y_1,z_1]$, where $$[[\underbrace{\widetilde{x_s},\widetilde{x_s}],...,\widetilde{x_s}}_{\hbox{i-times}}]=(1+a_{n-1})^{i-1}e_i+(*)e_{i+1}+ \dots+(*)e_{n-3} \hbox{ for } 3 \leq i \leq n-3,$$
obtaining the gradation: $V_{k_s}\oplus V_{2k_s}\oplus \dots \oplus V_{(n-3)k_s}\oplus V_{k_t}\oplus V_{k_u}\oplus V_{k_s+k_t}.$ Let us assume that the gradation has maximum length, therefore $k_s, k_t,k_u$ are pairwise different.
It is enough to consider the products $[z_2,y_1]$ and $[y_1,z_2]$ to prove that it is not possible that the gradation has maximum length.
Consider $[z_2,y_1]=[\widetilde{x_u},\widetilde{x_s}]=c_1[\widetilde{x_s},\widetilde{x_s}]=c_1y_2$. Since $[z_2,y_1] \in V_{k_s+k_u}$, $y_{2} \in V_{2k_s}$ and $k_s \neq k_u$, then we conclude $c_1=0$.

On the other hand $[y_1,z_2]=(1+c_1)e_2+(*)e_3 + \dots+ (*)e_{n-3}+c_{n-2}f_3=e_2+(*)e_3+\dots+(*)e_{n-3}+\beta_1f_3=Ay_2$ with $A \neq 0$. We also have $[y_1,z_2] \in V_{k_s+k_u}$ and $y_{2} \in V_{2k_s}$, therefore $k_u=k_s$, which is a contradiction with the assumption of maximum length. Hence there is no maximum length algebra in this subcase.

\

\fbox{Case 1.2:} If $[\widetilde{x_s},\widetilde{x_s}]$ and $[\widetilde{x_s},\widetilde{x_t}]$ are linearly dependent, i.e., $(b_1+b_{n-1})a_{n-2}=1+a_{n-1}.$ Note that we can assert that $a_{n-2} \neq 0$ and $b_1+b_{n-1} \neq 0$ from the assumption $1+a_{n-1} \neq 0$.

\

\textbf{Case 1.2.1:} If $[\widetilde{x_s},\widetilde{x_t}]$ and $[\widetilde{x_s},\widetilde{x_u}]$ are linearly independent, we distinguish two possibilities:

\

\emph{\textbf{If $c_{n-2} \neq 0,$}}
let us take the basis composed of the following vectors $y_1=\widetilde{x_s},$ $y_i=[y_{i-1},y_1]$ for $2 \leq i \leq n-3,$ $z_1=\widetilde{x_t},$ $z_2=\widetilde{x_u}$ and
$z_3=[\widetilde{x_s},\widetilde{x_u}]=[y_1,z_2]$, and the maximum length gradation
$V_{k_s}\oplus V_{2k_s}\oplus \dots \oplus V_{(n-3)k_s}\oplus V_{k_t}\oplus V_{k_u}\oplus V_{k_s+k_u}$. A contradiction will be obtained
by computing $[y_1,z_1].$ Since
\begin{align*}
 [y_1,z_1]&=\underbrace{(b_1+b_{n-1})}_{\neq 0}e_2+ (*)e_3+ \dots+(*)e_{n-3}+\underbrace{c_{n-2}}_{\neq 0}f_3 \in V_{k_s+k_t},\\
  y_2=[\widetilde{x}_s,\widetilde{x}_s]&=\underbrace{(1+a_{n-1})}_{\neq 0}e_2+(*)e_3+ \dots + (*)e_{n-3}+\underbrace{a_{n-2}}_{\neq 0}f_3 \in V_{2k_s},\\
 z_3=[\widetilde{x_s},\widetilde{x_u}]&=\underbrace{(1+c_1)}_{\neq 0}e_2+(*)e_3 + \dots+ (*)e_{n-3}+\underbrace{c_{n-2}}_{\neq 0}f_3 \in V_{k_s+k_u},
 \end{align*}
then either $[y_1,z_1]=Ay_2$ with $A \neq 0$ or $[y_1,z_1]=Bz_3$ with $B \neq 0$. By properties of the gradation, we achieve $k_s=k_t$ in the first case and $k_u=k_t$ in the other case. Both equalities contradict the maximum length of $\widetilde{L^1}.$

\

\emph{\textbf{If $c_{n-2}=0,$}} we can assume that $1+c_1 \neq 0.$ Otherwise we would have:
\begin{align*}
[\widetilde{x_s},\widetilde{x_s}]&=(1+a_{n-1})e_2+(*)e_3+ \dots + (*)e_{n-3}+a_{n-2}f_3,\\
[\widetilde{x_s},\widetilde{x_t}]&=(b_1+b_{n-1})e_2+(*)e_3 + \dots+ (*)e_{n-3}+f_3= \alpha[\widetilde{x_s},\widetilde{x_s}],\\
[\widetilde{x_s},\widetilde{x_u}]&=(*)e_3 + \dots+ (*)e_{n-3},
\end{align*}
and the others products would be linear combinations of these. Therefore, it would not be possible to generate the element $y_2$ or $z_3$ in the new basis. Hence, by taking the new basis $y_1=\widetilde{x_s},$
 $y_2=[\widetilde{x_s},\widetilde{x_u}],$ $y_i=[y_{i-1},y_1]$ for $3 \leq i \leq n-3,$ $z_1=\widetilde{x_t},$ $z_2=\widetilde{x_u}$ and
 $z_3=[\widetilde{x_s},\widetilde{x_t}]=[y_1,z_1]$, we get a contradiction as above, by calculating the product $[y_1,y_1]$.
 Since $$[y_1,y_1]=\underbrace{(1+a_{n-1})}_{\neq 0}e_2+(*)e_3+ \dots + (*)e_{n-3}+\underbrace{a_{n-2}}_{\neq 0}f_3,$$
 there are two possibilities, either
 $[y_1,y_1]=Ay_2$ with $A\neq 0$ or $[y_1,y_1]=Bz_3$ with $B\neq 0$. Analogously to the above case, both equalities contradict the hypothesis of maximum length of $\widetilde{L^1}.$

\

\textbf{Case 1.2.2:} If $[\widetilde{x_s},\widetilde{x_t}]$ and $[\widetilde{x_s},\widetilde{x_u}]$ are linearly dependent it is not possible to
construct a homogeneous basis, because all the products of the generators can be written as follows:
\begin{align*}
[\widetilde{x_t},\widetilde{x_t}]&=b_1[\widetilde{x_s},\widetilde{x_t}]=b_1\alpha[\widetilde{x_s},\widetilde{x_s}],\\
[\widetilde{x_u},\widetilde{x_u}]&=c_1[\widetilde{x_s},\widetilde{x_u}]=c_1\beta[\widetilde{x_s},\widetilde{x_s}],\\
[\widetilde{x_s},\widetilde{x_t}]&=\alpha[\widetilde{x_s},\widetilde{x_s}],\\
[\widetilde{x_t},\widetilde{x_s}]&=b_1[\widetilde{x_s},\widetilde{x_s}], \\
[\widetilde{x_s},\widetilde{x_u}]&=\beta[\widetilde{x_s},\widetilde{x_s}], \\
[\widetilde{x_u},\widetilde{x_s}]&= c_1[\widetilde{x_s},\widetilde{x_s}],\\
[\widetilde{x_t},\widetilde{x_u}]&=b_1[\widetilde{x_s},\widetilde{x_u}]=b_1\beta[\widetilde{x_s},\widetilde{x_s}], \\
[\widetilde{x_u},\widetilde{x_t}]&=c_1[\widetilde{x_s},\widetilde{x_t}]=c_1\alpha[\widetilde{x_s},\widetilde{x_s}],
\end{align*}
such that, all of them are linearly dependent of $[\widetilde{x_s},\widetilde{x_s}].$

\

$\bullet$ Case 2: If $1+a_{n-1}=0.$\\

\

\fbox{Case 2.1:} If $1+c_1 \neq 0$.\\

\

\textbf{Case 2.1.1:} If $[\widetilde{x_s},\widetilde{x_t}]$ and $[\widetilde{x_s},\widetilde{x_u}]$ are linearly independent, we take the new basis $y_1=\widetilde{x_s}$, $z_1=\widetilde{x_t}$, $z_2=\widetilde{x_u},$ $y_2=[\widetilde{x_s},\widetilde{x_u}]=[y_1,z_2],$ $y_i=[y_{i-1},z_2]$ for $3\leq i \leq n-3,$ and $z_3=[\widetilde{x_s},\widetilde{x_t}]=[y_1,z_1],$ where
$$y_{i+1}=[[\widetilde{x_s},\underbrace{\widetilde{x_u}],...,\widetilde{x_u}}_{i-times}]=(1+c_1)^i e_{i+1}+\dots+(*)e_{n-3}$$
with $2 \leq i \leq n-4,$ giving rise to the following gradation: $V_{k_s}\oplus V_{k_s+k_u}\oplus V_{k_s+2k_u}\oplus \dots \oplus V_{k_s+(n-4)k_u}\oplus V_{k_t}\oplus V_{k_u} \oplus V_{k_t+k_s}$ of maximum length.

In order to prove that there is no maximum length algebra in this subcase it is enough to study the values of $k_u,$ $k_t$ and $k_s$ such
as the gradation considered previously has maximum length. From properties of the gradation it is easy to check that the gradation is connected if and only if
$k_u=\pm 1.$ Without loss of generality we can assume $k_u=1.$ Let us study the values of $k_t$ and $k_s$.

\begin{itemize}

\item[$ $]    \textbf{Case a: If $k_s>0.$}
    \begin{figure}[ht]
\centering{
\subfigure[Subcase a.1]{\includegraphics[width=0.5\textwidth]{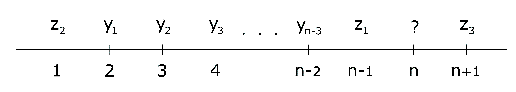}}
\hspace{0.05\textwidth}
\subfigure[Subcase a.2]{\includegraphics[width=0.4\textwidth]{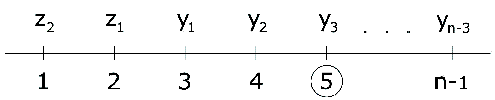}}
\hspace{0.05\textwidth}
}
\end{figure}
 \begin{figure}[ht]
\centering{
\subfigure[Subcase a.3]{\includegraphics[width=0.5\textwidth]{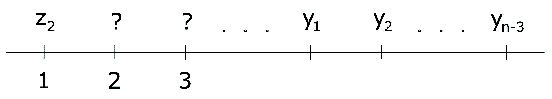}}
}
\end{figure}

     \textbf{Subcase a.1: If $k_s=2.$} Under these hypothesis we have $k_t=n-1,$ $V_{k_t+k_s}=V_{n+1}$ and $V_{n}=<0>.$ Hence the considered gradation is not connected.

      \textbf{Subcase a.2: If $k_s=3.$} Then $k_t=2$ and $V_{k_t+k_s}=V_5=<z_3,y_3>$, so the length of the gradation in not maximum.

   \textbf{Subcase a.3: If $k_s>3.$} Then there is some $V_p$ with $2 \leq p \leq k_s$ such that $V_p=<0>,$ which contradicts the connectedness of the gradation.

\item[$ $]   \textbf{Case b: If $k_s<0.$}
    \begin{figure}[ht]
\centering{
\subfigure[ Subcase b.1 and Subcase b.2]{\includegraphics[width=0.8\textwidth]{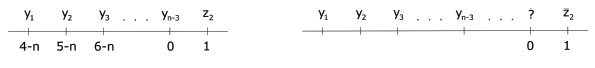}}
}
\end{figure}

     \textbf{Subcase b.1: $k_s=4-n$}. Under these hypothesis and by connectedness we have either $k_t=3-n$ or $k_t=2$. If  $k_t=3-n$ then
     $V_{k_t+k_s}=V_{7-2n}.$ Since $n\geq 7$, hence $V_{7-2n}$ is vanish and the gradation is not connected. On the other hand, if $k_t=2$ then
     $V_{k_t+k_s}=V_{6-n}=<y_3,z_3>$, which gives a contradiction with the assumption of maximum length.

      \textbf{Subcase b.2: $k_s \neq 4-n$}. This subcase never gives a maximum length gradation because the subspace $V_0$ is always vanish ($z_1\notin V_0$).

\end{itemize}

\textbf{Case 2.1.2}: If $[\widetilde{x_s},\widetilde{x_t}]$ and $[\widetilde{x_s},\widetilde{x_u}]$ are linearly dependent we have:
$$(1):\begin{cases}
[\widetilde{x_s},\widetilde{x_s}]=(*)e_3+ \dots+(*)e_{n-3}+a_{n-2}f_3\\
[\widetilde{x_s},\widetilde{x_t}]=(b_1+b_{n-1})e_2+(*)e_3+ \dots+(*)e_{n-3}+f_3=\alpha[\widetilde{x_s},\widetilde{x_u}],\\
[\widetilde{x_s},\widetilde{x_u}]=\underbrace{(1+c_1)}_{\neq 0}e_2+(*)e_3+ \dots+(*)e_{n-3}+c_{n-2}f_3. \\
\end{cases}$$
and the others products are linearly dependent of these. Therefore we can assume $a_{n-2} \neq 0,$ otherwise it would not be possible to get a basis because
$$\begin{array}{cc}
det \left(\begin{array}{cc}
b_1+b_{n-1}& 1 \\
1+c_1 & c_{n-2}
\end{array}\right)= 0\ \ \hbox{ and }&
det \left(\begin{array}{ccc}
1 & 0 & -1 \\
b_1 & 1 & b_{n-1}\\
c_1 & c_{n-2} & 1
\end{array}\right)= 0,
\end{array}$$
which implies that $\widetilde{x_s},$ $\widetilde{x_t}$ and $\widetilde{x_u}$ are linearly dependent.

 Let us take the new basis $y_1=\widetilde{x_s},$ $z_1=\widetilde{x_t},$ $z_2=\widetilde{x_u},$ $y_2=[\widetilde{x_s},\widetilde{x_u}]=[y_1,z_2],$
 $y_i=[y_{i-1},z_2]$ for $3 \leq i \leq n-3$ and $z_3=[y_1,y_1],$ where
$$y_{i}=[[\widetilde{x_s},\underbrace{\widetilde{x_u}],...,\widetilde{x_u}}_{(i-1)-times}]=(1+c_1)^{i-1} e_{i}+(*)e_{i+1}+\dots+(*)e_{n-3}$$
 with $3 \leq i \leq n-3.$ Its associated maximum length gradation is: $V_{k_s}\oplus V_{k_s+k_u}\oplus V_{k_s+2k_u}\oplus \dots \oplus V_{k_s+(n-4)k_u}\oplus V_{k_t}\oplus V_{k_u}
 \oplus V_{2k_s}$. Since $[\widetilde{x_s},\widetilde{x_t}]$ is linearly dependent
  of $[\widetilde{x_s},\widetilde{x_u}]$ and from properties of the gradation we conclude that $k_t=k_u,$ which is a contradiction with the assumption of
 maximum length.

\fbox{Case 2.2:} If $1+c_1=0,$ it is not possible to construct a basis because all the products of the generators can be written as follows:
\begin{align*}
[\widetilde{x_s},\widetilde{x_s}]&=(*)e_3+ \dots+(*)e_{n-3}+a_{n-2}f_3\\
[\widetilde{x_t},\widetilde{x_t}]&=b_1[\widetilde{x_s},\widetilde{x_t}]=(*)e_3+ \dots+(*)e_{n-3}+b_1a_{n-2}f_3,\\
[\widetilde{x_u},\widetilde{x_u}]&=c_1[\widetilde{x_s},\widetilde{x_u}]=(*)e_3+ \dots+(*)e_{n-3}+c_1c_{n-2}f_3,\\
[\widetilde{x_s},\widetilde{x_t}]&=(*)e_3+ \dots+(*)e_{n-3}+f_3,\\
[\widetilde{x_t},\widetilde{x_s}]&=b_1[\widetilde{x_s},\widetilde{x_s}]=(*)e_3+ \dots+(*)e_{n-3}+b_1a_{n-2}f_3, \\
[\widetilde{x_s},\widetilde{x_u}]&=(*)e_3+ \dots+(*)e_{n-3}+c_{n-2}f_3, \\
[\widetilde{x_u},\widetilde{x_s}]&= c_1[\widetilde{x_s},\widetilde{x_s}]=(*)e_3+ \dots+(*)e_{n-3}+c_1a_{n-2}f_3,\\
[\widetilde{x_t},\widetilde{x_u}]&=b_1[\widetilde{x_s},\widetilde{x_u}]=(*)e_3+ \dots+(*)e_{n-3}+b_1c_{n-2}f_3, \\
[\widetilde{x_u},\widetilde{x_t}]&=c_1[\widetilde{x_s},\widetilde{x_t}]=(*)e_3+ \dots+(*)e_{n-3}+c_1f_3,\\
\end{align*}
such that, the element $y_2$ can not be generated in this case. Therefore the proof is closed.
\end{dem}

\subsection{Split case}

\

\

This section is devoted to the study of the 3-filiform  Leibniz algebras of maximum length, whose naturally graded algebras are split. Furthermore we will focus our attention in the non standard families. The definitions of standard and non standard algebras are the following:

\begin{defn}
Let $\ll$ be a split maximum length algebra and let $k$ be an integer where $\ll=N_1 \oplus N_2 \oplus \dots \oplus N_k$. The algebra $\ll$ is called standard if $N_1$, $N_2$, $\dots$ and $N_k$ are algebras of maximum length. Otherwise the algebra $\ll$ is called non standard.
\end{defn}

\begin{exam}
The list of standard 3-filiform Leibniz algebras of maximum length consists of the following algebras: ma\-xi\-mum length null-filiform Leibniz algebras $\oplus \CC^3$, maximum length filiform Leibniz algebras $\oplus \CC^2$ and maximum length 2-filiform Leibniz algebras $\oplus \CC$. Note that these algebras have already been studied in \cite{Ayupov} and \cite{J.Lie.Theory2}.
\end{exam}

 Due to the previous example, we reduce our study to the non standard families, i.e., we study the extension of the naturally graded filiform non-Lie Leibniz algebras $\oplus \CC^2$ and the naturally graded 2-filiform non-Lie Leibniz algebras $\oplus \CC$. It should be remarked that the null-filiform case will not be studied because its extension always gives a standard algebra. The Lie case has already been done in \cite{3-filiform}, where the classification is presented in the following theorem:
\begin{thm}
Let $\ll$ be  a $(n+1)$-dimensional non standard 3-filiform Leibniz algebra whose associated naturally graded algebra is a Lie algebra. Then $n$ is odd and the algebra $\ll$ is isomorphic to the maximum length Lie algebra:
$$N: \begin{cases}
[e_{i-1},e_0]=e_i, \quad 2 \leq i \leq n-2,\\
[e_{n-3},e_1]=-e_{n-1},\\
[e_{n-4},e_2]=e_{n-1},\\
[e_{i},e_{n-2-i}]=(-1)^{i-1}e_{n-1}, \quad 3 \leq i \leq \lfloor\frac{n-3}{2}\rfloor,\\
[f_1,e_0]=e_{n-1}.
\end{cases}$$

\end{thm}

\fbox{2-Filiform case}

\

Cabezas, Camacho and Rodr\'{\i}guez gave the classification of the naturally graded 2-filiform non-Lie Leibniz algebras in \cite{J.Lie.Theory2}. They proved that, up to isomorphisms, there are two algebras under these hypothesis, which are not split. These algebras are defined by the following table of multiplications:
\begin{align*}
KF_4:&\begin{cases}
[e_i,e_1]=e_{i+1}, \quad 1 \leq i \leq n-3,\\
[e_1,e_{n-1}]=e_n+\alpha_3e_3+ \dots + \alpha_{n-2}e_{n-2},\\
[e_{n-1},e_{n-1}]=\beta_3e_3+\beta_4e_4+ \dots + \beta_{n-2}e_{n-2},\\
[e_i,e_{n-1}]=\beta_{i,i+2}e_{i+2}+\beta_{i,i+3}e_{i+3}+ \dots + \beta_{i,n-2}e_{n-2}, \quad 2 \leq i \leq n-4,\\
[e_n,e_{n-1}]=\gamma_4e_4+ \dots + \gamma_{n-2}e_{n-2}.\\
\end{cases}\\[2mm]
KF_5:&\begin{cases}
[e_i,e_1]=e_{i+1}, \quad 1 \leq i \leq n-3,\\
[e_1,e_{n-1}]=e_2+e_n+\alpha_3e_3+ \dots + \alpha_{n-2}e_{n-2},\\
[e_{n-1},e_{n-1}]=\beta_3e_3+\beta_4e_4+ \dots + \beta_{n-2}e_{n-2},\\
[e_i,e_{n-1}]=e_{i+1}+\beta_{i,i+2}e_{i+2}+\beta_{i,i+3}e_{i+3}+ \dots + \beta_{i,n-2}e_{n-2}, \quad 2 \leq i \leq n-4,\\
[e_n,e_{n-1}]=\gamma_4e_4+ \dots + \gamma_{n-2}e_{n-2}.\\
\end{cases}
\end{align*}
Due to the above classification we obtain the following theorem:
\begin{thm}
Let $\ll$ be a $(n+1)-dimensional$ 3-filiform non-Lie Leibniz algebra of maximum length, whose associated naturally graded algebra is $KF_4\oplus \mathbb{C}.$ Then $\ll$ is isomorphic to either $M$ or one of the algebra of the family $M^{1,\alpha}$:
$$\begin{array}{cc}
M: \begin{cases}
[y_i,y_1]=y_{i+1}, \quad 1 \leq i \leq n-3,\\
[y_1,y_{n-1}]=y_n,\\
[z_1,y_{n-1}]=y_{n-2},
\end{cases}&
M^{1,\alpha}:\begin{cases}
[y_i,y_1]=y_{i+1}, \quad 1 \leq i \leq n-3,\\
[y_1,y_{n-1}]=y_n,\\
[y_{n-1},z_{1}]=y_{n-2},\\
[z_1,y_{n-1}]=\alpha y_{n-2}, \quad \alpha\in \mathbb{C}.\\
\end{cases}
\end{array}$$
 \end{thm}

 \begin{dem}
The extension of the algebra $KF_4\oplus \mathbb{C}$, via the natural gradation, is:
$$\widetilde{\ll}:\begin{cases}
[e_i,e_1]=e_{i+1}+(*)e_{i+2}+\dots+(*)e_{n-2}, &1\leq i \leq n-3,\\
[e_1,e_{n-1}]=e_n+(*)e_3+\dots+(*)e_{n-2},\\
[e_{n-1},e_{n-1}]=(*)e_{3}\dots+(*)e_{n-2},\\
[e_i,e_{n-1}]=(*)e_{i+2}+\dots+(*)e_{n-2}, &2\leq i\leq n-4,\\
[e_n,e_{n-1}]=(*)e_4+\dots+(*)e_{n-2},\\
[e_i,f_1]=(*)e_{i+2}+\dots +(*)e_{n-2}, &1 \leq i \leq n-4,\\
[e_{n-1},f_1]=(*)e_3+\dots+(*)e_{n-2},\\
[e_n,f_1]=(*)e_4+ \dots + (*)e_{n-2},\\
[f_1,e_i]=(*)e_{i+2}+\dots+(*)e_{n-2}, &1 \leq i \leq n-4, \\
[f_1,e_{n-1}]=(*)e_3+\dots+(*)e_{n-2},\\
[f_1,e_n]=(*)e_4+ \dots + (*)e_{n-2},\\
[f_1,f_1]=(*)e_{3}+ \dots+ (*)e_{n-2},
\end{cases}$$
where the asterisks $(*)$ denote the corresponding coefficients in the products.
We are going to get the homogenous basis by considering the generators:
\begin{align*}
\widetilde{x_s}&=e_1+ \sum_{i=2}^n a_ie_i+b_1f_1,\\
\widetilde{x_t}&=e_{n-1}+ \sum_{i=1, i\neq n-1}^n A_ie_i+B_1f_1,\\
\widetilde{x_u}&=f_1+ \sum_{i=1}^n \alpha_ie_i.
\end{align*}

Let us consider the following products, since they will be very useful in the rest of the proof:
\begin{align*}
[\widetilde{x_s},\widetilde{x_s}]&=e_2+(*)e_3+\dots+ (*)e_{n-2}+a_{n-1}e_n,\\
[\widetilde{x_s},\widetilde{x_t}]&=A_1e_2+(*)e_3+ \dots + (*)e_{n-2}+e_n,\\
[[\underbrace{\widetilde{x_s},\widetilde{x_s}], \dots,\widetilde{x_s}}_{\hbox{i-times}}]&=e_{i}+(*)e_{i+1}+ \dots+(*)e_{n-2}, \hbox{ with } 3 \leq i \leq n-2.
\end{align*}

Let us take the homogeneous basis $y_1=\widetilde{x_s},$ $y_i=[y_{i-1},y_1]$ for $2 \leq i \leq n-2,$ $y_{n-1}=\widetilde{x_t},$
$y_n=[y_1,y_{n-1}],$ $z_1=\widetilde{x_u}$ and the associated maximum length gradation $V_{k_s}\oplus V_{2k_s} \oplus \dots \oplus
 V_{(n-2)k_s} \oplus V_{k_t} \oplus V_{k_t+k_s} \oplus V_{k_u}.$ This gradation is connected if and only if $k_s=\pm 1.$ Without loss of generality
  we can assume $k_s=1$ (the case $k_s=-1$ is analogous). We are going to continue the proof studying the possible values that the subindices
   $k_t$ and $k_u$ can get  to obtain a maximum length gradation.

\

$\bullet$ Case 1: If $k_t>0$, there are the following possibilities:

    \begin{figure}[ht]
\centering{
\subfigure[Subcase 1.1]{\includegraphics[width=0.3\textwidth]{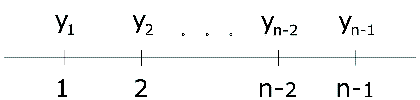}}
\hspace{0.05\textwidth}
\subfigure[Subcase 1.2]{\includegraphics[width=0.4\textwidth]{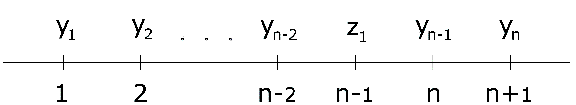}}
\hspace{0.05\textwidth}
}
\end{figure}

 \begin{figure}[ht]
\centering{
\subfigure[Subcase 1.3]{\includegraphics[width=0.5\textwidth]{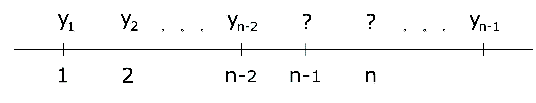}}
}
\end{figure}

\textbf{Subcase 1.1: If $k_t =n-1,$} then from the connectedness of the gradation we obtain $k_u=0$ or $k_u=n+1.$ But $k_u=0$ is not possible because
$z_1 \in V_{k_u}$ and $z_1$ is a generator. If $k_u=n+1$ then $[y_i,z_1] \in V_{n+1+i}=<0>$, $[z_1,y_i] \in V_{n+1+i}=<0>$ for $1 \leq i \leq n-2$. Moreover
$[z_1,y_{n-1}]$ and $[y_{n-1},z_1]$ belong to the subspace $V_{2n}=<0>$. On the other hand $[z_1,y_n],$ $[y_n,z_1] \in V_{2n+1}=<0>$ and
$[z_1,z_1] \in V_{2n+2}=<0>.$ Then $z_1 \in Cent(\widetilde{\ll}),$ giving rise to a standard algebra.

\textbf{Subcase 1.2: If $k_t=n.$} By a similar previous reason, it is clear that $[y_i,z_1]=[z_1,y_i]=0$ for $3 \leq i \leq n$ because
those products belong to $V_{n-1+i}=<0>.$ Moreover $[z_1,y_1],[y_1,z_1] \in V_n=<y_{n-1}>$, but this is not possible since $y_{n-1}$ is a generator,
such that $y_{n-1} \in \widetilde{\ll}\setminus \widetilde{\ll}^2$, while $[z_1,y_1],[y_1,z_1] \in \widetilde{\ll}^2.$

Since $y_2 \in R(\widetilde{\ll})$ and from the Leibniz identity we affirm that $[z_1,y_2]=[y_2,z_1]=0.$ Finally, since $[z_1,z_1] \in V_{2n-2}=<0>$, we conclude that $z_1 \in Cent(\widetilde{\ll}),$ such that, the obtained algebra is standard.

\textbf{Subcase 1.3: If $k_t>n,$} the gradation is not connected because either $V_{n-1}=<0>$ or $V_{n}=<0>$.\\

$\bullet$ Case 2: If $k_t<0$, we distinguish:
\begin{figure}[ht]
\centering{
\subfigure[Subcase 2.1]{\includegraphics[width=0.4\textwidth]{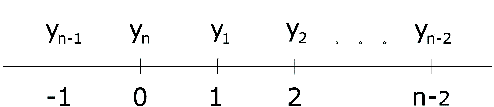}}
\hspace{0.05\textwidth}
\subfigure[Subcase 2.2]{\includegraphics[width=0.5\textwidth]{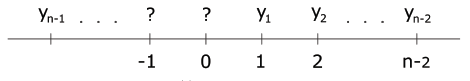}}
\hspace{0.05\textwidth}
}
\end{figure}

\textbf{Subcase 2.1: If $k_t=-1,$} then from the connectedness of the gradation either $k_u=-2$ or $k_u=n-1$. In the first case we are going to prove that the obtained algebra is standard because $z_1 \in Cent(\widetilde{\ll})$. From properties of the gradation we have $[y_i,z_1], [z_1,y_i] \in V_{i-2}=<y_{i-2}>$ for $2 \leq i \leq n-2,$ but from the law of $\widetilde{\ll}$ we know $[y_i,z_1]=\alpha_1e_{i+1}+(*)e_{i+2}+ \dots+(*)e_{n-2}$ and $[z_1,y_i]=(*)e_{i+2}+\dots+(*)e_{n-2}.$ Therefore $[y_i,z_1]=0$ and $[z_1,y_i]=0$ are concluded for $2 \leq i \leq n-2.$ On the other hand $[y_1,z_1], [z_1,y_1] \in V_{-1}=<y_{n-1}>,$ which is not possible because $y_{n-1}$ is a generator of $\widetilde{\ll}.$ Finally, since $[z_1,y_{n-1}], [y_{n-1},z_1] \in V_{-3}=<0>$ and  $y_n \in Cent(\widetilde{\ll}),$ we have finally proved that $z_1 \in Cent(\widetilde{\ll})$, such that, the algebra is standard.
\

If $k_u=n-1,$ $[z_1,y_i]=[y_i,z_1]=0$ for $1 \leq i \leq n-2$ because they belong to $V_{n-1+i}=<0>$ and $[z_1,z_1]=0$ because $[z_1,z_1]\in V_{2n-2}=<0>.$
 Moreover since $y_n \in Cent(\widetilde{\ll})$ we assert that $[z_1,y_n]=[y_n,z_1]=0.$ On the other hand, from properties of the gradation we can write $[z_1,y_{n-1}]=\alpha y_{n-2}$ and
$[y_{n-1},z_1]=\beta y_{n-2}.$ Due to $\{y_2,y_3, \dots, y_{n-2}\} \in R(\widetilde{\ll})$, $y_n \in Cent(\widetilde{\ll})$ and the above calculations, it is enough to
compute the products: $[y_{n-2},y_1],$ $[y_{n-1},y_1]$ and $[y_{i},y_{n-1}],$ for $2 \leq i \leq n-1,$ to know the law of $\widetilde{\ll}$ in the homogeneous basis.

From properties of the gradation it is clear that $[y_{n-2},y_1] \in V_{n-1}=<z_1>.$ But from the definition of the descending central sequence we have $[y_{n-2},y_1] \in \widetilde{\ll}^2$ and $z_1 \in \widetilde{\ll}\setminus \widetilde{\ll}^2$. Hence we get $[y_{n-2},y_1]=0.$ By the same arguments, it can be concluded that $[y_1,y_{n-2}]=0.$ Furthermore $[y_{n-1},y_{n-1}]=0$ since $[y_{n-1},y_{n-1}] \in V_{-2}=<0>.$

Besides it can be proved that $[y_3,y_{n-1}] \in V_2=<y_2>$ and from the law of $\widetilde{\ll}$ we can write $[y_3,y_{n-1}]=A_1e_4+(*)e_5+ \dots+e_{n-2}=A_1y_4.$ So we conclude $A_1=0$ and $[y_3,y_{n-1}]=0,$ because if $A_1 \neq 0,$ then  $V_2\supseteq [y_3,y_{n-1}]=A_1y_4 \in V_4$, which contradicts the assumption of maximum length of the gradation. Analogously we get $[y_{n-1},y_3]=0.$ In addiction, as
$[y_{n-1},y_1]=A_1e_2+(*)e_3+\dots+(*)e_{n-2}+A_1a_{n-1}e_n= (*)e_3+\dots+(*)e_{n-2}$ and $[y_{n-1},y_1] \in V_{0}=<y_n>,$ then holds $[y_{n-1},y_1]=0.$

In summary, the obtained law of the maximum length algebra is:
$$\widetilde{\ll}:\begin{cases}
[y_i,y_1]=y_{i+1}, \quad 1 \leq i \leq n-3,\\
[y_1,y_{n-1}]=y_n,\\
[y_i,y_{n-1}]=\gamma_iy_{i-1}, \quad 2 \leq i \leq n-2,\\
[z_1,y_{n-1}]=\alpha y_{n-2},\\
[y_{n-1},z_1]=\beta y_{n-2}.
\end{cases}$$

Finally, by using the program of the Leibniz identity it is easy to prove that $\gamma_i=0$ for $2 \leq i \leq n-2.$ Further by considering the dimension of $R(\ll),$ we assume $\beta=0$ or $\beta=1$. On the one hand if $\beta=0$ it is necessary that $\alpha \neq 0$ and by a trivial change of basis we can take $\alpha=1$. This gives rise to $M$. On the other hand ($\beta=1$), by using the program of the isomorphism we obtain the family $M^{1,\alpha}$, with $\alpha \in \mathbb{C}.$

\

\textbf{Subcase 2.2: If $k_t \neq -1,$} we only attain standard algebras or not connected gradations, by similar arguments as in previous cases.
 \end{dem}

 \begin{thm}
Let $\ll$ be a $(n+1)-dimensional$ 3-filiform non-Lie Leibniz algebra, whose associated naturally graded algebra is $KF_5\oplus \mathbb{C}.$ Then $l(\ll)\leq n.$
 \end{thm}
 \begin{dem}
 The proof is achieved by using a similar reasoning to that followed in the previous theorem: to take a homogeneous basis and the
 associated maximum gradation, to use the properties of the gradation and the above programs.
 \end{dem}

\fbox{Filiform case}

\

 Ayupov and Omirov in \cite{Ayupov} obtained the classification of naturally graded filiform non-Lie Leibniz algebras in arbitrary dimension. They proved that, up to isomorphisms, there are three algebras for each dimension $n$. We are going to show only one, because the other ones are either a split algebra or a Lie algebra.
 $$NGF_1: \begin{cases}
 [e_1,e_1]=e_3,\\
 [e_i,e_1]=e_{i+1}, \quad 2 \leq i \leq n-1.
 \end{cases}$$

 Extending the algebra $NGF_1 \oplus \mathbb{C}^2$, via the natural gradation, the following result is attained:
 \begin{thm}
Let $\ll$ be a $(n+2)$-dimensional 3-filiform non-Lie Leibniz algebra of maximum length, whose associated naturally graded algebra is
 $NGF_1\oplus \mathbb{C}^2,$ with $n \geq 8.$ Then $l(\ll)\leq n+1.$

 \end{thm}

 \begin{dem}
As in previous proofs, the first step is to consider the extension of the algebra $NGF_1 \oplus \mathbb{C}^2$, by using its natural gradation, and to get a
homogeneous basis derived from the generators
\begin{align*}
\widetilde{x_s}&=e_1+ \sum_{i=2}^n a_ie_i+b_1f_1+b_2f_2,\\
\widetilde{x_t}&=e_{2}+ \sum_{i=1, i\neq 2}^{n}
A_ie_i+B_1f_1+B_2f_2,\\
\widetilde{x_u}&=f_1+ \sum_{i=1}^n \alpha_ie_i+\beta_2f_2,\\
\widetilde{x_v}&=f_2+ \sum_{i=1}^n \gamma_ie_i+\mu_1f_1.
\end{align*}

The main products of these generators are:
$$(2): \begin{cases}
 [\widetilde{x_s},\widetilde{x_s}]&=(1+a_2)e_3+ (*)e_4 + \dots+ (*)e_{n-2},\\
 [\widetilde{x_t},\widetilde{x_s}]&=(1+A_1)e_3+ (*)e_4 + \dots+ (*)e_{n-2},\\
 [\widetilde{x_u},\widetilde{x_s}]&=(\alpha_1+\alpha_2)e_3+ (*)e_4 + \dots+ (*)e_{n-2},\\
 [\widetilde{x_v},\widetilde{x_s}]&=(\gamma_1+\gamma_2)e_3+ (*)e_4 + \dots+ (*)e_{n-2},\\
 \end{cases}$$
because the other products are linearly dependent of these.

 The next step is to assume that the associated gradation with that basis has maximum length. Let us see in details.

\

 $\bullet$ Case 1: If $1+a_2 \neq 0,$ we take the basis $y_1=\widetilde{x_s},$ $y_2=\widetilde{x_t},$ $y_3=[y_1,y_1],$ $y_i=[y_{i-1},y_1]$ for
 $4 \leq i \leq n$, $z_1=\widetilde{x_u}$ and $z_2=\widetilde{x_v}$ and the associated gradation $V_{k_s}\oplus V_{2k_s} \oplus
  V_{(n-1)k_s} \oplus V_{k_t} \oplus V_{k_u} \oplus V_{k_v},$ whose length is maximum. Note that $$(3):[\underbrace{[\widetilde{x_s},\widetilde{x_s}], \dots,\widetilde{x_s}}_{\hbox{i-times}}]=(1+a_2)e_{i+1}+ (*)e_{i+2}+ \dots+ (*)e_{n-2}, \ 2 \leq i \leq n-1.$$

  From (2) and (3) we conclude that $[y_2,y_1]$ is linearly dependent of $y_3.$ In addiction
  from properties of the gradation $[y_2,y_1] \in V_{k_t+k_s}$ and $y_3 \in V_{2k_s}.$ Finally, by the hypothesis of maximum length we know that $k_s \neq k_t$. These facts imply $[y_2,y_1]=0$, hence $A_1=-1$ (see (2)). On the other hand $[y_1,y_2]=A_1(1+a_2)e_3+ (*)e_4+ \dots+ (*)e_{n-2}=-y_3,$ $[y_1,y_2] \in V_{k_s+k_t}$ and $y_3 \in V_{2k_s}$, then $V_{k_s+k_t}=V_{2k_s}$, such that, $k_s=k_t$
   which is not possible. We conclude that there is no maximum length algebra in this case.

\

 $\bullet$ Case 2: If $1+a_2=0,$ we have to distinguish the following cases:

\

 \textbf{Subcase 2.1: If} $A_1=0,$ we take the new basis $y_1=x_s,$ $y_2=x_t,$ $y_i=[y_{i-1},y_1]$ for $3 \leq i \leq n,$ $z_1=x_u$ and  $z_2=x_v.$ The associated maximum length gradation is: $V_{k_s}\oplus V_{k_t} \oplus V_{k_t+k_s}\oplus V_{k_t+2k_s} \oplus \dots \oplus V_{k_t+(n-2)k_s} \oplus V_{k_u} \oplus V_{k_v}.$

  We now consider all the possible product in the new basis, obtaining the following law:
 $$\begin{cases}
[y_1,y_1]=y_3,\\
[y_i,y_1]=y_{i+1}, \qquad 3 \leq i \leq n-1,\\
[y_i,y_{2}]=P_iy_{i+4}, \quad 2 \leq i \leq n-4,\\
[y_{i},z_{1}]=Q_iy_{i+2}, \quad 2 \leq i \leq n-2,\\
[y_i,z_{2}]=R_iy_{i+3}, \quad 2 \leq i \leq n-3.
\end{cases}$$

It is clear to see, by induction on $i$ and by calculating the Leibniz identity on
$[[y_i,y_1],y_2],$ $[[y_i,y_1],z_1]$ and $[[y_i,y_1],z_2]$, that $P_i=Q_i=R_i=0$ for $i\geq 3,$ respectively. Moreover by applying the program of the Leibniz identity, we prove that
 $A_2=B_2=D_2=0$ (for more details see the fo\-llo\-wing Web site: http://personal.us.es/jrgomez). Therefore the obtained algebra is standard.\\

\textbf{Subcase 2.2:} If $A_1 \neq 0 \wedge A_1 \neq -1,$ then we can take the same previous homogeneous basis. We get a contradiction with the assumption of maximum length by considering the product $[y_2,y_2].$ From the law of $\widetilde{\ll}$ we have $[y_2,y_2]=A_1(A_1+1)e_3 +(*)e_4+\dots +(*)e_{n-2}=A_1y_3.$ Since we had assumed that $A_1 \neq 0$, $[y_2,y_2] \in V_{2k_t}$ and $y_3 \in V_{k_t+2k_s}$, then it can be concluded that $k_s=k_t$, which is a contradiction. \\

 \textbf{Subcase 2.3:} If $A_1 \neq 0 \wedge A_1=-1,$ since $\widetilde{x_u}$ and $\widetilde{x_v}$ play a symmetric role, we can assume
 that $\alpha_1+\alpha_2 \neq 0.$ Otherwise it was not possible to construct a homogeneous basis generated by
 $\widetilde{x_u}$, $\widetilde{x_t},$  $\widetilde{x_u}$ and $\widetilde{x_v}$. Therefore we take the basis $y_1=x_s,$ $y_2=x_u,$ $y_i=[y_{i-1},y_1]$ for $3 \leq i \leq n,$ $z_1=x_t$ and $z_2=x_v$.

 If $\alpha_1 \neq 0,$ then $[y_2,y_2]=[\widetilde{x_u},\widetilde{x_u}]=\alpha_1(\alpha_1+\alpha_2)e_3+(*)e_4+ \dots+(*)e_{n-2}=\tau y_3 \neq 0.$ Since $[y_2,y_2]\in V_{2k_u},$ $y_3 \in V_{k_t+k_s}$ and $\tau \neq 0$, $k_t=k_s$ is achieved, which is not possible because it contradicts the maximum length. Therefore $\alpha_1=0$, $\alpha_2 \neq 0$ and $[\widetilde{x_u},\widetilde{x_t}]=-\alpha_2e_3+(*)e_4+ \dots+(*)e_{n-2},$ obtaining, by a similar way, the same contradiction.
 \end{dem}

\section{Applications of maximum length.}

As we mentioned in the introduction, the algebras of maximum length allow to study some cohomological properties easily, such as the space of derivations and the first cohomology group (see \cite{loday-piras}).

We have centred on the computational support again in order to tackle these cohomological pro\-per\-ties. We will use a third program, the program of derivations, that allows to determinate a basis of the space of derivations of an algebra of maximum length. From here, the cohomology study can be easily completed by using similar arguments as in \cite{Ayupov}, \cite{Dz1}--\cite{Reyes1}, \cite{Omirov2}, \cite{Ve}. As in the other programs, the implementation is presented in low and fixed dimension (see \cite{mitesis} for more details). Then we will formulate the generalizations, proving by induction the results for arbitrary fixed dimension.

\begin{pr}

\

\begin{itemize}
\item $dim(\mathcal{D}er(N))=3\displaystyle\frac{n-1}{2}+7.$\\

\item $dim(\mathcal{D}er(M))=n+6.$\\

\item $dim(\mathcal{D}er(M^{1,\alpha}))=n+5.$

\end{itemize}
\end{pr}

\begin{dem}
The proof is carried out by using the program of derivations, whose calculations are presented with a step-by-step explanation in the following Web site:
http://personal.us.es/jrgomez.
\end{dem}

\begin{cor}
\

\begin{itemize}
\item $dim(\mathcal{H}^1(N,N))=\displaystyle\frac{n+19}{2}$.\\

\item $dim(\mathcal{H}^1(M,M))=n+4$.\\

\item $dim(\mathcal{H}^1(M^{1,\alpha},M^{1,\alpha}))=n+2$.
\end{itemize}
\end{cor}

\begin{dem}
The proof is carried out by using the characterization $H^{1}(\ll,\ll)=\mathcal{D}er(\ll)\setminus \mathcal{I}nn(\ll),$ where $\mathcal{I}nn(\ll)$ denotes the set of the inner derivations of $\ll.$
\end{dem}

\

\

{\sc Luisa M. Camacho, Elisa M. Ca\~{n}ete, Jos\'{e} R. G\'{o}mez.}  Dpto. Matem\'{a}tica Aplicada I.
Universidad de Sevilla. Avda. Reina Mercedes, s/n. 41012 Sevilla.
(Spain), e-mail: \emph{lcamacho@us.es}, \emph{elisacamol@us.es}, \emph{jrgomez@us.es}

\

{\sc Bakhrom A. Omirov.} Institute of Mathematics. National University of Uzbekistan,
F. Hodjaev str. 29, 100125, Tashkent
(Uzbekistan), e-mail: \emph{omirovb@mail.ru}


\begin{thebibliography}{12}


\bibitem{Ayupov}{\rm Sh.A. Ayupov and B.A. Omirov,} {\rm On some classes of nilpotent Leibniz algebras}, (Russian) \emph{Sibirsk. Math.} \textbf{42 (1)}
(2001) 18--29; translation in \emph{Siberian Math. J.} \textbf{42 (1)} (2001) 15--24.

\bibitem{J.Lie.Theory2}{\rm J.M. Cabezas, L.M. Camacho and I.M. Rodr\'{\i}guez,} {\rm On filiform and 2-filiform Leibniz algebras of maximum length}, \emph{Journal of Lie Theory} \textbf{18}
(2008) 335--350.
\bibitem{Pastor}{\rm J.M. Cabezas and E. Pastor,} {\rm Naturally graded $p$-filiform Lie algebras in arbitrary finite dimension}, \emph{Journal of Lie Theory} \textbf{15}
(2005) 379--391.



\bibitem{3-filiform}{\rm L.M. Camacho, E.M. Ca\~{n}ete, J.R. G\'{o}mez and B.A. Omirov,} {\rm 3-filiform Leibniz algebras of maximum length, whose naturally graded algebras are Lie algebras}, \emph{Linear and Multilinear Algebra} \textbf{59(9)} (2011) 1039-1059.

    \bibitem{JSC} {\rm L.M. Camacho, J.R. G\'{o}mez, A.J. Gonz\'{a}lez and B.A. Omirov,} {\rm Naturally graded quasi-filiform Leibniz algebras},
\emph{Journal of Symbolic Computation} \textbf{44(5)} (2009) 527-539.


\bibitem{NGp-F} {\rm L.M. Camacho, J.R. G\'{o}mez, A.J. Gonz\'{a}lez and B.A. Omirov,} {\rm The classification of naturally graded $p$-filiform Leibniz algebras},
\emph{Communications in Algebra} \textbf{39 (1)} (2011) 153--163.

\bibitem{Casas} {\rm J.M. Casas and M. Ladra,} {\rm Non-abelian tensor product of Leibniz algebras and an exact sequence in Leibniz homology,} \emph{Comm. Algebra} \textbf{31 (9)} (2003) 4639--4646.

\bibitem{mitesis}{\rm E.M. Ca\~{n}ete,} {\rm Algebras de Leibniz de longitud maxima,} PhD Thesis, Universidad de Sevilla, 2012. (http://www.educacion.es/teseo).

\bibitem{Dz1}
{\rm A.S. Dzhumadil'daev,} {\it Cohomologies of colour Leibniz algebras: pre-simplicial approach.} Lie theory and its applications in physics, III (Clausthal, 1999), 124--136, World Sci. Publ., River Edge, NJ, 2000.


\bibitem{Cohomology}
{\rm  B.L. Feigin and  D.B. Fuks,} {\rm Cohomology of Lie groups and Lie algebras, II,} 125--223, Encyclopadedia Math. Sci., 21, Springer, Berlin, 2000.

\bibitem{Reyes1}{\rm J.R. G\'{o}mez, A. Jim\'{e}nez-Merch\'{a}n and J. Reyes,} {\rm Maximum length filiform Lie algebras}, \emph{Extracta Math.} \textbf{16(3)} (2001) 405--421.



%



\bibitem{loday1} {\rm J.L. Loday,} {\it Cyclic homology,}  Springer-Verlag, Berlin, 1992.


\bibitem{loday-piras} {\rm J.L. Loday and T. Pirashvili,} {\rm Universal enveloping algebras of Leibniz algebras and (co)homology,}
\emph{Math. Ann.} \textbf{296(1)} (1993) 139-158.

\bibitem{Omirov2}
{B.A. Omirov,} {\rm On derivations of filiform Leibniz algebras}, \emph{Math. Notes} \textbf{77 (5)} (2005) 733--742.

\bibitem{Ve}
{\rm M. Vergne,} {\rm Cohomologie des alg\`{e}bres de Lie
nilpotentes. Application`a l'\'{e}tude de la variet\'{e} des
alg\`{e}bres de Lie nilpotentes}, \emph{Bull. Soc. Math. France}
\textbf{98} (1970) 81--116.
\end{thebibliography}
\end{document}